\documentclass[11pt]{article}

\usepackage{latexsym}
\usepackage{amssymb}
\usepackage{amsthm}
\usepackage{amscd}
\usepackage{amsmath}
\usepackage[all]{xy}
\usepackage{graphicx}

\newtheorem{theorem}{Theorem}[section]

\newtheorem{lemma}{Lemma}[section]

\newtheorem{proposition}{Proposition}[section]
\newtheorem{corollary}{Corollary}[section]

\newtheorem{conjecture}{Conjecture}[section]

\xyoption{dvips}

\numberwithin{equation}{section}

\newcommand{\FF}{\mathbb{F}}

\newcommand{\RR}{\mathbb{R}}
\newcommand{\CC}{\mathbb{C}}

\def\GL{\mathrm{GL}}

\def\SL{\mathrm{SL}}

\def\Irr{\mathrm{Irr}}

\newcommand{\bZ}{\boldsymbol{Z}}

\newcommand{\U}{\mathrm{U}}
\newcommand{\gO}{\mathrm{O}}
\newcommand{\Sp}{\mathrm{Sp}}
\newcommand{\GSp}{\mathrm{GSp}}
\newcommand{\GO}{\mathrm{GO}}
\newcommand{\CSp}{\mathrm{CSp}}

\newcommand{\SO}{\mathrm{SO}}
\newcommand{\ep}{\varepsilon}

\newcommand{\bG}{\boldsymbol{G}}

\makeatletter
\renewcommand{\@makefnmark}{\mbox{\textsuperscript{}}}
\makeatother

\def\adots{\mathinner{\mkern2mu\raise0pt\hbox{.}  
\mkern2mu\raise4pt\hbox{.}\mkern1mu
\raise7pt\vbox{\kern7pt\hbox{.}}\mkern1mu}}

\allowdisplaybreaks[1]

\begin{document}

\bibliographystyle{amsplain}

\title{Character degree sums and real representations of finite classical groups of odd characteristic}
\author {C. Ryan Vinroot}
\date{}

\maketitle

\begin{abstract}
Let $\FF_q$ be a finite field with $q$ elements, where $q$ is the power of an odd prime, and let $\GSp(2n, \FF_q)$ and $\GO^{\pm}(2n, \FF_q)$ denote the symplectic and orthogonal groups of similitudes over $\FF_q$, respectively.  We prove that every real-valued irreducible character of $\GSp(2n, \FF_q)$ or $\GO^{\pm}(2n, \FF_q)$ is the character of a real representation, and we find the sum of the dimensions of the real representations of each of these groups.  We also show that if $\bG$ is a classical connected group defined over $\FF_q$ with connected center, with dimension $d$ and rank $r$, then the sum of the degrees of the irreducible characters of $\bG(\FF_q)$ is bounded above by $(q+1)^{(d+r)/2}$.  Finally, we show that if $\bG$ is any connected reductive group defined over $\FF_q$, for any $q$, the sum of the degrees of the irreducible characters of $\bG(\FF_q)$ is bounded below by $q^{(d-r)/2}(q-1)^r$.  We conjecture that this sum can always be bounded above by $q^{(d-r)/2}(q+1)^r$.
\\
\\
2000 {\it Mathematics Subject Classification:}  20C33, 20G40
\end{abstract}




\section{Introduction}

Given a finite group $G$, and an irreducible complex representation $(\pi, V)$ of $G$ which is self-dual (that is, has a real-valued character), one may ask whether $(\pi, V)$ is a real representation.  Frobenius and Schur gave a method of answering this question by introducing an invariant, which we denote $\ep(\pi)$, or $\ep(\chi)$ when $\chi$ is the character of $\pi$ (called the Frobenius-Schur indicator), which gives the value $0$ when $\chi$ is not real-valued, $1$ when $\pi$ is a real representation, and $-1$ when $\chi$ is real-valued but $\pi$ is not a real representation (see \cite[Chapter 4]{Is76}).  For example, if $S_n$ is the symmetric group on $n$ elements, it is known that all of the representations of $S_n$ are real (and, in fact, rational).  It follows from results of Frobenius and Schur that the sum of all of the dimensions of the irreducible representations of $S_n$ is equal to the number of elements in $S_n$ which square to the identity element.

Now consider the group $\GL(n, \FF_q)$ of invertible linear transformations on an $n$-dimensional vector space over the field $\FF_q$ with $q$ elements, where $q$ is the power of a prime $p$.  Not every irreducible character of $\GL(n, \FF_q)$ is real-valued, however, every real-valued irreducible character of $\GL(n, \FF_q)$ is the character of a real representation.  That is, $\ep(\pi) = 0$ or $1$ for every irreducible representation $\pi$ of $\GL(n, \FF_q)$.  This follows from Ohmori's result \cite{Oh77} that every irreducible character of $\GL(n, \FF_q)$ has rational Schur index 1 (also proved by Zelevinsky \cite{Ze81}), and a direct proof is given by Prasad \cite{Pr98}.  This result also follows from a theorem of Gow \cite{Go83}, which states that when $q$ is the power of an odd prime, the sum of the dimensions of the irreducible representations of $\GL(n, \FF_q)$ is equal to the number of symmetric matrices in the group (also obtained by Klyachko \cite{Kl83} and Macdonald \cite{Mac95} for all $q$).  Gow's proof actually implies that the {\em twisted} Frobenius-Schur indicator (see Section \ref{realreps}) of an irreducible representation of $\GL(n, \FF_q)$, with respect to the transpose-inverse automorphism, is always $1$.  Using the result that every real-valued character of $\GL(n, \FF_q)$ is the character of a real representation, in Section \ref{realsums}, we compute the sum of the dimensions of the real representations of the group.  This gives one way to compare the size of the set of real representations to the size of the set of all representations of $\GL(n, \FF_q)$.

In Section \ref{realreps} of this paper, we consider the groups of symplectic and orthogonal similitudes over a finite field, denoted $\GSp(2n, \FF_q)$ and $\GO^{\pm}(2n, \FF_q)$ respectively, where $q$ is the power of an odd prime.  The main results of Section \ref{realreps}, which are Theorems \ref{GSpReal} and \ref{GOReal}, state that, as in the case of the group $\GL(n, \FF_q)$, every real-valued irreducible character of $\GSp(2n, \FF_q)$ or $\GO^{\pm}(2n, \FF_q)$ is the character of a real representation.  There is a very direct proof in the case of the groups $\GO^{\pm}(2n, \FF_q)$, which comes from the fact that every irreducible representation of the orthogonal groups $\gO^{\pm}(2n, \FF_q)$, where $q$ is the power of an odd prime, is a real representation, another result of Gow \cite{Go85}.  Also in \cite{Go85}, Gow proves that a real-valued character of the symplectic group $\Sp(2n, \FF_q)$, where $q$ is the power of an odd prime, is the character of real representation if and only if the central element $-I$ acts trivially, and in particular, has irreducible real-valued characters which are not the characters of real representations.  It may seem to be a surprising result, then, that there are no such characters for the group of symplectic similitudes $\GSp(2n, \FF_q)$ by Theorem \ref{GSpReal}.  However, the situation is similar for the finite special linear group, $\SL(n, \FF_q)$, where in the case that $n$ is congruent to $2$ mod $4$ and $q$ is congruent to $1$ mod $4$ the group has real-valued characters which are not characters of real representations (see \cite{Go81}), while all irreducible real-valued characters of $\GL(n, \FF_q)$ are characters of real representations.  So, in this instance, the relationship between the characters of $\GSp(2n, \FF_q)$ and $\Sp(2n, \FF_q)$ is similar to the relationship between the characters of $\GL(n, \FF_q)$ and $\SL(n, \FF_q)$.

In Section \ref{realsums}, we apply the results of Section \ref{realreps} to compute expressions for the sum of the dimensions of the real representations of the groups $\GSp(2n, \FF_q)$ and $\GO^{\pm}(2n, \FF_q)$, in Theorems \ref{GSpRealSum} and \ref{GORealSum}.  By results of Frobenius and Schur and Theorems \ref{GSpReal} and \ref{GOReal}, we need only count the number of elements in each group which square to the identity element.  In the resulting expressions, we obtain a term which is the sum of all of the dimensions of the representations of $\Sp(2n, \FF_q)$ and $\gO^{\pm}(2n, \FF_q)$.  This sum has been computed for $\Sp(2n, \FF_q)$ by Gow in the case that $q$ is congruent to $1$ mod $4$ in \cite{Go85}, and by the author in the case that $q$ is congruent to $3$ mod $4$ in \cite{Vi05}.  This sum has not been computed in the case of the orthogonal groups $\gO^{\pm}(n, \FF_q)$, however, and so we do this in Section \ref{OrthogSums}, with the result given in Theorem \ref{OrthSum}.  We also give several corollaries for the special orthogonal groups in Section \ref{OrthogSums}.

In the last three sections, we focus on the following result of Kowalski \cite[Proposition 5.5]{Ko08}, which he obtains in the context of sieving applications.

\begin{theorem} [Kowalski] \label{KoThm} Let $\bG$ be a split connected reductive group with connected center over $\bar{\FF}_q$, defined over $\FF_q$.  Let $d$ be the dimension of $\bG$, let $r$ be the rank of $\bG$, let $W$ be the Weyl group of $\bG$, and let $G = \bG(\FF_q)$.  Then, the sum of the degrees of the irreducible complex characters of $G$ is bounded above as follows:
$$ \sum_{\chi \in {\rm Irr}(G)} \chi(1) \leq (q+1)^{(d+r)/2} \left(1 + \frac{2r|W|}{q-1} \right).$$
\end{theorem}

Kowalski also notes that in the cases of the groups $\GL(n, \FF_q)$ and $\GSp(2n, \FF_q)$ ($q$ the power of an odd prime), by results of Gow \cite{Go83} and the author \cite{Vi05}, the factor $1 + \frac{2r|W|}{q-1}$ may be removed from the bound in Theorem \ref{KoThm}.  In Section \ref{IneqLemmas}, we prove several inequalities which we apply to improve Theorem \ref{KoThm} in the cases of orthogonal groups.  In particular, in Theorem \ref{MainBound}, which is the main result of Section \ref{Bound}, we prove that the factor $1 + \frac{2r|W|}{q-1}$ may be removed from the bound in Theorem \ref{KoThm} for any connected classical group (removing the restriction of being split) with connected center which is defined over $\FF_q$, with $q$ the power of an odd prime.  In the case of the unitary group $\U(n, \FF_{q^2})$, this follows from a specific formula for the sum of the degrees of the irreducible characters due to Thiem and the author \cite{ThVi07}.  The other groups to check are the special orthogonal groups $\SO(2n+1, \FF_q)$ and the connected orthogonal similitude groups $\GO^{\pm, \circ}(2n, \FF_q)$, which is where our results from Sections \ref{OrthogSums} and \ref{IneqLemmas} are applied.  

Finally, in Section \ref{Lower}, we turn to the general case of a finite reductive group with connected center.  We use the Gelfand-Graev character to prove, in Proposition \ref{LowerProp}, that if $\bG$ is a connected reductive group with connected center defined over $\FF_q$ (for any $q$) with dimension $d$ and rank $r$, then the sum of the degrees of the irreducible complex characters of $\bG(\FF_q)$ is bounded below by $q^{(d-r)/2}(q-1)^r$.  We notice that in the cases of the finite general linear, unitary, and symplectic similitude groups (for $q$ odd), the sum of the degrees of the irreducible characters can actually be bounded above by $q^{(d-r)/2}(q+1)^r$, and we conjecture that this is the case for all finite reductive groups with connected center.
\\
\\
{\bf Acknowledgments. }  The author thanks Julio Brau for checking some of the calculations made in Section \ref{OrthogSums}.

\section{Real representations of finite classical groups} \label{realreps}

Let $G$ be a finite group and $\sigma$ an automorphism of $G$ such that $\sigma^2$ is the identity.  Let $(\pi, W)$ be an irreducible complex representation of $G$, and suppose that ${^\sigma \pi} \cong \hat{\pi}$, where $\hat{\pi}$ is the contragredient representation of $\pi$ and ${^\sigma \pi}$ is defined by ${^\sigma \pi}(g) = \pi(\sigma(g))$.  If $\chi$ is the character of $\pi$, note that ${^\sigma \pi} \cong \hat{\pi}$ if and only if ${^\sigma \chi} = \bar{\chi}$.  From the isomorphism ${^\sigma \pi} \cong \hat{\pi}$, there must exist a nondegenerate bilinear form $B_{\sigma}: W \times W \rightarrow \CC$, unique up to scalar by Schur's Lemma, such that
$$ B_{\sigma}({^\sigma \pi}(g) u, \pi(g) v) = B_{\sigma}(u, v), \;\; \text{ for all } g \in G, \, u, v \in W.$$
Since switching the variables of $B_{\sigma}$ gives a bilinear form with the same property, and the bilinear form is unique up to scalar, then $B_{\sigma}$ must be either symmetric or skew-symmetric.  That is, we have
$$B_{\sigma}(u, v) = \varepsilon_{\sigma}(\pi) B_{\sigma}(v, u) \;\; \text{ for all } u, v \in W,$$
where $\ep_{\sigma}(\pi) = \pm 1$.  If ${^\sigma \pi} \not\cong \hat{\pi}$, then define $\ep_{\sigma}(\pi) = 0$.  We will also write $\ep_{\sigma}(\pi) = \ep_{\sigma}(\chi)$.  In the case that $\sigma$ is the identity automorphism, then $\ep_{\sigma}(\pi) = \ep(\pi)$ is just the classical Frobenius-Schur indicator of $\pi$.  In this case, if $\chi$ is the character of $\pi$, then $\ep(\chi) = 0$ exactly when $\chi$ is not real-valued, and when $\chi = \bar{\chi}$, then $\ep(\chi) = 1$ if $\chi$ is the character of a real representation, and $\ep(\chi) = -1$ otherwise.  In the case that $\sigma$ is an order $2$ automorphism, the invariant $\ep_{\sigma}(\pi)$ is called the {\em twisted} Frobenius-Schur indicator of $\pi$, which was first considered by Mackey \cite{Ma58}, and later studied in more detail by Kawanaka and Matsuyama \cite{KaMa90}.  For any finite group $G$, let $\Irr(G)$ denote the set of complex irreducible characters of $G$.  The twisted Frobenius-Schur indicators satisfy the following identity (see \cite{KaMa90}), which reduces to the classical Frobenius-Schur involution formula in the case that $\sigma$ is trivial:
$$ \sum_{\chi \in \Irr(G)} \ep_{\sigma}(\chi) \chi(1) = \big| \{ g \in G \, \mid \, \sigma(g) = g^{-1} \} \big|.$$

If $G$ is a finite group with $(\pi, W)$ an irreducible complex representation with character $\chi$, and $z$ is a central element of $G$, then by Schur's Lemma $\pi(z)$ acts by a scalar on $W$.  We let $\omega_{\pi}(z)$, or $\omega_{\chi}(z)$, denote this scalar.  If $\sigma$ is an order $2$ automorphism of $G$ and $z$ is a central element of $G$, we define $G(\sigma, z)$ to be the following central extension of $G$:
$$ G(\sigma, z) = \langle G, \tau \; \mid \; \tau^2 = z, \tau^{-1} g \tau = \sigma(g) \text{ for all } g \in G \rangle.$$
Note that $G$ is an index $2$ subgroup of $G(\sigma, z)$.  In the case that $z = 1$, we just have that $G(\sigma, 1)$ is the split extension of $G$ by $\sigma$.  We have the following properties of the representations of $G$ and $G(\sigma, z)$, and their Frobenius-Schur indicators, which were proven in \cite[Section 2]{Vi05}.

\begin{proposition} \label{induce} Let $\chi$ be an irreducible character of $G$, and $\chi^{+}$ the induced character ${\rm Ind}_G^{G(\sigma, z)}(\chi)$.  Then:
\begin{enumerate}
\item $\chi^{+}$ is irreducible if and only if ${^\sigma \chi} \neq \chi$.  In this case, we have $\ep(\chi^{+}) = \ep(\chi) + \omega_{\chi}(z) \ep_{\sigma}(\chi)$.
\item If ${^\sigma \chi} = \chi$, then $\chi^+ = \psi_1 + \psi_2$, where $\psi_1$ and $\psi_2$ are irreducible characters of $G(\sigma, z)$.  In this case, we have $\ep(\psi_1) + \ep(\psi_2) = \ep(\chi) + \omega_{\chi}(z) \ep_{\sigma}(\chi)$.  Also, each $\psi_i$ is an extension of the character $\chi$, and $\ep(\psi_1) = \ep(\psi_2)$.
\end{enumerate}
\end{proposition}

In statement 2 in the proposition above, although the last claim is not specifically proven in \cite[Section 2]{Vi05}, it follows from the definition of the induced character that each $\psi_i$ is an extension of $\chi$, and it follows from the same calculation as in \cite[Lemma 2.3]{Vi05} that we have $\ep(\psi_i) = \frac{1}{2}(\ep(\chi) + \omega_{\chi}(z)\ep_{\sigma}(\chi))$ for $i=1,2$, so that $\ep(\psi_1) = \ep(\psi_2)$.  The next lemma is used in the proof of the first main result of this section.

\begin{lemma} \label{MainLemma}
Let $G$ be a finite group, and $\sigma$ an order 2 automorphism of $G$.  Suppose that every character of $G(\sigma, z)$ is real-valued, and that every irreducible character $\chi$ of $G$ satisfies $\ep_{\sigma}(\chi) = \omega_{\chi}(z)$.  Then, every real-valued character $\chi$ of $G$ satisfies $\ep(\chi) = 1$, and every irreducible character $\psi$ of $G(\sigma, z)$ satisfies $\ep(\psi) = 1$.
\end{lemma}
\begin{proof} Let $\chi$ be an irreducible real-valued character of $G$, so that $\ep(\chi) = \pm 1$.  We are assuming that $\ep_{\sigma}(\chi) = \omega_{\chi}(z) = \pm 1$, so that in particular we have ${^\sigma \chi} = \bar{\chi}$.  Since $\chi$ is real-valued, then ${^\sigma \chi} = \chi$, and by Proposition \ref{induce}, $\chi^+ = \psi_1 + \psi_2$, where
$$ \ep(\psi_1) + \ep(\psi_2) = \ep(\chi) + \omega_{\chi}(z) \ep_{\sigma}(\chi).$$
Now, we have $\omega_{\chi}(z) \ep_{\sigma}(\chi) = 1$ by assumption, and $\ep(\chi) = \pm 1$, $\ep(\psi_1) = \ep(\psi_2) = \pm 1$, by Proposition \ref{induce} and by assumption.  The only possibility is $\ep(\chi) = 1$ and $\ep(\psi_i) = 1$.

For the second statement, every irreducible character of $G(\sigma, z)$ is either extended or induced from an irreducible character of $G$, since $G$ is an index $2$ subgroup.  We have just checked that whenever $\psi$ is an irreducible of $G(\sigma, z)$ which is extended from $G$, then $\ep(\psi) = 1$.  Now assume that $\psi = {\rm Ind}_G^{G(\sigma, z)}(\chi)$ for some irreducible $\chi$ of $G$.  By Proposition \ref{induce}, we must have ${^\sigma \chi} \neq \chi$, and so $\chi \neq \bar{\chi}$ since we are assuming $\ep_{\sigma}(\chi) = \pm 1$.  From Proposition \ref{induce}(1), we have $\ep(\psi) = \ep(\chi) + \omega_{\chi}(z) \ep_{\sigma}(\chi)$.  Now, $\ep(\chi) = 0$ since $\chi \neq \bar{\chi}$, and $\omega_{\chi}(z) \ep_{\sigma}(\chi) = 1$ by assumption.  Thus $\ep(\psi) = 1$.
\end{proof}

Let $\FF_q$ be a finite field with $q$ elements, where $q$ is the power of an odd prime.  Let $V = \FF_q^{2n}$, and let $\langle \cdot, \cdot \rangle$ be a nondegenerate skew-symmetric form on $V$ (of which there is only one equivalence class, by \cite[Theorem 2.10]{Gr02}).  The group of all elements $g$ of $\GL(V)$ which leave $\langle \cdot, \cdot \rangle$ invariant up to a scalar multiple is called the {\em symplectic group of similitudes on $\FF_q^{2n}$} (or the {\em conformal symplectic group on $\FF_q^{2n}$}) which we will denote as $\GSp(2n, \FF_q)$ (this group is also denoted as $\CSp(2n, \FF_q)$).  That is, $\GSp(2n, \FF_q)$ is the group of elements $g$ of $\GL(2n, \FF_q)$ such that $\langle gv, gw \rangle = \mu(g) \langle v, w \rangle$ for all $v, w \in V$, where $\mu(g) \in \FF_q^{\times}$ is a scalar depending only on $g$.  Then $\mu: \GSp(2n, \FF_q) \rightarrow \FF_q^{\times}$ is a character called the {\em similitude character}, and the symplectic group over $\FF_q$, denoted $\Sp(2n, \FF_q)$, is the kernel of $\mu$.

Now consider the case that $G = \GSp(2n, \FF_q)$, where $q$ is the power of an odd prime, and we let $\sigma$ be the order $2$ automorphism of $G$ defined by $\sigma(g) = \mu(g)^{-1} g$, where $\mu$ is the similitude character.  We have the following results for the group $\GSp(2n, \FF_q)$ proven previously by the author, the first statement in \cite[Theorem 6.2]{Vi05}, and the second in \cite{Vi04}.

\begin{proposition} \label{factor} 
Let $G = \GSp(2n, \FF_q)$, where $q$ is the power of an odd prime.  Then:
\begin{enumerate}
\item Every irreducible character $\chi$ of $G$ satisfies $\ep_{\sigma}(\chi) = \omega_{\chi}(-I)$, where $\sigma(g) = \mu(g)^{-1} g$.
\item For any $g \in G$, we may write $g = h_1 h_2$ such that $h_1^2 = I$, $\mu(h_1) = -1$, and $h_2^2 = \mu(g) I, \mu(h_2) = - \mu(g)$.
\end{enumerate}
\end{proposition}

We are now ready to prove the following result. 

\begin{theorem} \label{GSpReal} Let $G = \GSp(2n, \FF_q)$, and $q$ the power of an odd prime.  Then:
\begin{enumerate}
\item[i.] If $\chi$ is an irreducible real-valued character of $G$, then $\ep(\chi) = 1$.  

\item[ii.] If $\sigma$ is the order $2$ automorphism of $G$ defined by $\sigma(g) = \mu(g)^{-1} g$, then every irreducible character $\psi$ of $G(\sigma, -I)$ satisfies $\ep(\psi) = 1$.
\end{enumerate}
\end{theorem}
\begin{proof} Since $\ep_{\sigma}(\chi) = \omega_{\chi}(-I)$ for every irreducible character $\chi$ of $G$, then by Lemma \ref{MainLemma}, it is enough to show that every irreducible character of $G(\sigma, -I)$ is real-valued.  Equivalently, we must show that every element of $G(\sigma, -I)$ is conjugate to its inverse.  We prove this by applying the factorization given by Proposition \ref{factor}(2).

First, let $g \in G \subset G(\sigma, -I)$.  Write $g = h_1 h_2$ as in Proposition \ref{factor}(2).  Then, 
$$g^{-1} = h_2^{-1} h_1^{-1} = \mu(g)^{-1} h_2 h_1.$$
If we conjugate $g$ by $h_1 \tau$, we obtain
$$ (h_1 \tau) g (h_1 \tau)^{-1} = h_1 (\tau g \tau^{-1}) h_1 = h_1 (\mu(g)^{-1} h_1 h_2) h_1 = \mu(g)^{-1} h_2 h_1 = g^{-1},$$
and so $g$ is conjugate to $g^{-1}$.  Now let $g \tau \in G(\sigma, -I) \setminus G$, and again write $g = h_1 h_2$ as in Proposition \ref{factor}(2).  In this case,
$$(g\tau)^{-1} = \tau^{-1}\mu(g)^{-1} h_2 h_1 = -\mu(g) \tau h_2 h_1 = -\mu(g) \mu(h_2)^{-1} h_2 h_1 \tau = h_2 h_1 \tau,$$
since $\mu(h_2) = - \mu(g)$.  Now conjugate $g\tau$ by $h_1$, and we have
$$ h_1(g\tau)h_1 = h_2 \tau h_1 = h_2 h_1 \tau = (g \tau)^{-1}.$$
So, $g\tau$ is conjugate to its inverse, and all characters of $G(\sigma, -I)$ are real-valued.
\end{proof}

Gow \cite{Go83} proved that the split extension of $\GL(n, \FF_q)$ by the transpose-inverse automorphism has the property that all of its characters are characters of real representations.  This is analogous to our result for the central extension of $\GSp(2n, \FF_q)$ in Theorem \ref{GSpReal}(ii).

Now consider a nondegenerate symmetric form $\langle \cdot, \cdot \rangle$ on a vector space $V = \FF_q^n$ over a finite field $\FF_q$, where $q$ is odd.  Similar to the symplectic case, the {\em orthogonal group of similitudes on $\langle \cdot, \cdot \rangle$} ({\em or conformal orthogonal group}) is the group of elements $g \in GL(V)$ which leave the form $\langle \cdot, \cdot \rangle$ invariant up to a scalar multiple.  When the form does not need to be emphasized, we will denote this group by $\GO(n, \FF_q)$ (the group is sometimes denoted ${\rm CO}(n, \FF_q)$).  As in the symplectic case, we let $\mu$ denote the similitude character, where if $g \in \GO(n, \FF_q)$, then $\langle gu, gv \rangle = \mu(g) \langle u, v \rangle$ for all $u, v \in V$.  In particular, the orthogonal group $\gO(n, \FF_q)$ for the form $\langle \cdot, \cdot \rangle$ is the kernel of the homomorphism $\mu: \GO(n, \FF_q) \rightarrow \FF_q^{\times}$.

When $n = 2m+1$ is odd, any symmetric form on $V$ gives an isomorphic group (see \cite[Chapter 9]{Gr02}), and in fact the group $\GO(2m+1, \FF_q)$ is just the direct product of its center with the special orthogonal group $\SO(2m+1, \FF_q)$ \cite[Lemma 1.3]{Sh80}.  So, we will not consider the odd case.  When $n = 2m$ is even, however, there are two different equivalence classes of symmetric forms, split and nonsplit (see \cite[Chapter 9]{Gr02} or \cite[Section 15.3]{DiMi91}), giving two non-isomorphic groups denoted $\GO^+(2m, \FF_q)$ and $\GO^-(2m, \FF_q)$, respectively.  When considering both cases, and the type of form makes a difference in the result, we will write $\GO^{\pm}(2m, \FF_q)$.  Similarly, we let ${\rm O}^{+}(n, \FF_q)$ denote a split orthogonal group, ${\rm O}^-(n, \FF_q)$ the non-split orthogonal group, and in the case that $n$ is odd, we just write $\gO(n, \FF_q)$ for the unique orthogonal group.  When considering the split and non-split orthogonal groups at the same time, we use the notation $\gO^{\pm}(n, \FF_q)$.  We also use $\gO(n, \FF_q)$ for any of these orthogonal groups when the general case is considered.

We now study the real-valued characters of the groups $\GO^{\pm}(2n, \FF_q)$.  We could proceed in a fashion similar to the case of the group $\GSp(2n, \FF_q)$, by defining an automorphism $\sigma(g) = \mu(g)^{-1} g$ on $\GO^{\pm}(2n, \FF_q)$.  The exact same argument will work, since in \cite{Vi06} the author proved the appropriate factorization result, and $\ep_{\sigma}(\chi) = 1$ for every irreducible character $\chi$ of any orthogonal similitude group.  However, a more direct argument is possible here, coming from the result of Gow \cite[Theorem 1]{Go85} that for odd $q$, every irreducible character of any orthogonal group $\gO(n, \FF_q)$ has Frobenius-Schur indicator $1$.  We give this direct argument in the following, and the proof is very similar to that of \cite[Theorem 2]{Vi06}.

\begin{theorem} \label{GOReal} Let $G = \GO^{\pm}(2n, \FF_q)$, where $q$ is the power of an odd prime.  Then every irreducible real-valued character $\chi$ of $G$ satisfies $\ep(\chi) = 1$.
\end{theorem}
\begin{proof}  Let $\chi$ be an irreducible real-valued character of $G$, which is the character of the representation $(\pi, W)$.  Then there is a nondegenerate bilinear form $B$ on $W$, unique up to scalar, such that
$$B(\pi(g) u, \pi(g) v) = B(u, v) \quad \text{for all } g \in G, \, u, v \in W,$$
which must be either symmetric or skew-symmetric, and $\ep(\chi) = 1$ is equivalent to $B$ being a symmetric form.

Now let $Z$ be the center of $G$, which consists of all nonzero scalar matrices, where we have $\mu(bI) = b^2$ for any $b \in \FF_q^{\times}$.  Let $H = Z \cdot \gO^{\pm}(2n, \FF_q)$, which is an index $2$ subgroup of $G$, since $G/H \cong \FF_q^{\times}/ (\FF_q^{\times})^2$.  Note that any irreducible representation of $H$ is an extension of an irreducible of $\gO^{\pm}(2n, \FF_q)$, since $Z$ is central.

Since $H$ is an index $2$ subgroup of $G$, then $\pi$ is either extended or induced from an irreducible representation of $H$.  First suppose that $\pi$ is extended from an irreducible $\rho$ of $H$, with character $\psi$.  Then $\psi$ is the extension of a real-valued irreducible character $\theta$ of $\gO^{\pm}(2n, \FF_q)$, and say $\theta$ is the character of the representation $(\phi, W)$.  By the result of Gow \cite{Go85}, we have $\ep(\theta) = 1$, which means there is a nondegenerate bilinear form on $W$, unique up to scalar, which is invariant under the action of $G$ under $\phi$, which must be symmetric.  But the form $B$ is such a form, since $(\pi, W)$ is an extension of $(\rho, W)$.  Therefore, $B$ must be symmetric and $\ep(\chi) = 1$ in this case.

Now suppose that $\pi$ is induced from an irreducible representation of $H$, which means that $\pi$ restricted to $H$ is isomorphic to the direct sum of two irreducible representations, $(\rho_1, W_1)$ and $(\rho_2, W_2)$ of $H$, which are extended from irreducible representations $(\phi_1, W_1)$ and $(\phi_2, W_2)$ of $\gO^{\pm}(2n, \FF_q)$, respectively.  If $\theta_1$ is the character of $\phi_1$, then again by Gow's result, $\ep(\theta_1) = 1$, and there is a nondegenerate bilinear form on $W_1$, unique up to scalar, which is invariant under the action of $\phi_1$ and is symmetric.  We can make the same conclusion as in the previous case if we can show that the form $B$ is nondegenerate on $W_1$, since then $B$ would be symmetric on a nonzero subspace, and thus symmetric everywhere.  If $B$ is nondegenerate on $W_1 \times W_2$, then for $u \in W_1$, $v \in W_2$, and $g \in G$, we would have $B(\pi(g) u, \pi(g) v) = B(\phi_1(g) u_1, \phi_2(g) u_2)$, which would imply $\phi_1 \cong \hat{\phi_2} \cong \phi_2$.  This, in turn, would imply that $\rho_1 \cong \rho_2$, but it is impossible for an irreducible representation to restrict to the direct sum of $2$ isomorphic representations of an index $2$ subgroup \cite[Corollary 6.19]{Is76}.  Then $B$ must be degenerate on $W_1 \times W_2$, and so must be nondegenerate on $W_1 \times W_1$, since it is nondegenerate on $W$.  Thus, $B$ must be symmetric, and $\ep(\chi) = 1$.
\end{proof}

In summary, when $q$ is odd, all real-valued irreducible characters of the groups $\GL(n, \FF_q)$ (see the Introduction), $\GSp(2n, \FF_q)$, and $\GO^{\pm}(2n, \FF_q)$ have Frobenius-Schur indicator $1$, and the same is also true for the special orthogonal groups $\SO(2n+1, \FF_q)$ \cite[Theorem 2]{Go85}.  These are all examples of groups of $\FF_q$-points of classical groups with connected center (although the algebraic groups $\GO^{\pm}(2n)$ are disconnected).  Another example one might consider is the finite unitary group $\U(n, \FF_{q^2})$, but this is known to have irreducible characters with Frobenius-Schur indicator equal to $-1$ (see \cite{Oh96}, for example).

\section{Degree sums of real-valued characters} \label{realsums}

As mentioned in the Introduction, every real-valued irreducible character of the group $\GL(n, \FF_q)$ is the character of a real representation.  It follows from the classical Frobenius-Schur involution formula (see \cite[Chapter 4]{Is76}) that when $G = \GL(n, \FF_q)$, we have
$$ \sum_{\chi \in {\rm Irr}(G) \atop{ \chi \; \RR \text{-valued}}} \chi(1) = | \{ g \in G \, \mid \; g^2 = I \}|.$$
We may count the number of elements in $\GL(n, \FF_q)$ which square to the identity by summing the indices of the centralizer of an element in each order 2 conjugacy class.  We will make such counts for several groups, and the following notation will be helpful for the expressions obtained.  For any $q > 1$, and any integers $m, k \geq 0$ such that $m \geq k$, the {\em $q$-binomial} or {\em Gaussian binomial coefficients} are defined as
$$ \binom{m}{k}_q = \frac{(q^m - 1) (q^{m-1} - 1) \cdots (q^{m-k+1} - 1)}{(q^k -1)(q^{k-1} -1) \cdots (q-1)}.$$
It follows from induction and the identity $\binom{m}{k}_q = \binom{m-1}{k}_q + q^{m-k} \binom{m-1}{k-1}_q$, $m \geq 1$, that the Gaussian binomial coefficients are polynomials in $q$.  The Gaussian binomial coefficients have analogous properties to the standard binomial coefficients, for example, we have $\binom{m}{k}_q = \binom{m}{m-k}_q$.

We now count the number of elements in $\GL(n, \FF_q)$ which square to the identity, when $q$ is the power of an odd prime.  Each such element in $\GL(n, \FF_q)$ must have elementary divisors only of the form $x \pm 1$, and so is conjugate to a diagonal element with only $1$'s and $-1$'s on the diagonal.  An element in the conjugacy class of elements with exactly $k$ eigenvalues which are equal to $1$, and $n-k$ eigenvalues equal to $-1$, has centralizer isomorphic to $\GL(k, \FF_q) \times \GL(n-k, \FF_q)$.  Since $|\GL(m, \FF_q)| = q^{m(m-1)/2} \prod_{i=1}^m (q^i -1)$, then the size of the conjugacy class containing the elements whose square is $I$ and which have exactly $k$ eigenvalues equal to $1$ is
$$ \frac{|\GL(n, \FF_q)|}{|\GL(k, \FF_q) \times \GL(n-k, \FF_q)|} = q^{k(n-k)} \frac{\prod_{i=1}^n (q^i -1)}{\prod_{i=1}^k (q^i -1) \prod_{i=1}^{n-k} (q^i -1)} = q^{k(n-k)} \binom{n}{k}_q.$$
Summing the size of each conjugacy class, for $0 \leq k \leq n$, gives the result
\begin{equation} \label{GLReal}
\sum_{\chi \in {\rm Irr}(G) \atop{ \chi \; \RR \text{-valued}}} \chi(1) = \sum_{k=0}^n q^{k(n-k)} \binom{n}{k}_q.
\end{equation}
Using the fact that as a polynomial in $q$, the degree of $\binom{n}{k}_q$ is $k(n-k)$, we calculate that the sum in (\ref{GLReal}), as a polynomial in $q$, has degree $n^2/2$ when $n$ is even, and degree $(n^2 - 1)/2$ when $n$ is odd.  On the other hand, from a result of Gow \cite[Theorem 4]{Go83}, the sum of the degrees of all of the irreducible characters of $\GL(n,\FF_q)$ as a polynomial in $q$ has degree $(n^2 + n)/2$.  The ratio of the sum of the degrees of real-valued characters over the sum of all degrees in this case is roughly $q^{n/2}$ when $n$ is even and $q^{(n+1)/2}$ when $n$ is odd.  This gives a quick comparison in measure of the set of real-valued irreducible characters of $\GL(n,\FF_q)$ to the set of all irreducible characters, relative to degree size.

By Theorem \ref{GSpReal}, when $q$ is odd, every real-valued character of $\GSp(2n, \FF_q)$ is the character of a real representation, and so the sum of the degrees of the real-valued characters of $\GSp(2n, \FF_q)$ is equal to the number of elements in $\GSp(2n, \FF_q)$ which square to $I$.  Similar to the case of $\GL(n, \FF_q)$, we may count these elements to obtain an expression in $q$ for this character degree sum.

\begin{theorem} \label{GSpRealSum}
Let $q$ be the power of an odd prime, let $G = \GSp(2n, \FF_q)$, and let $H = \Sp(2n, \FF_q)$.  Then, the sum of the degrees of the real-valued characters of $G$ is given by
\begin{align*}
\sum_{\chi \in {\rm Irr}(G) \atop{ \chi \; \RR \text{-valued}}} \chi(1)  = \sum_{k=0}^n q^{2k(n-k)} \binom{n}{k}_{q^2} + \sum_{ \chi \in {\rm Irr}(H)} \chi(1) & = \sum_{k=0}^n q^{2k(n-k)} \binom{n}{k}_{q^2} + \frac{|\Sp(2n, \FF_q)|}{|\GL(n, \FF_q)|} \\
 & = \sum_{k=0}^n q^{2k(n-k)} \binom{n}{k}_{q^2} + q^{n(n+1)/2} \prod_{i=1}^n (q^i + 1).  
\end{align*}
\end{theorem}
\begin{proof} We must count the number of elements of $\GSp(2n, \FF_q)$ such that $g^2 = I$.  Such an element must satisfy $\mu(g) = \pm 1$.  If $\mu(g) = 1$, then we have $g \in \Sp(2n, \FF_q)$, and we may count such elements using results on the conjugacy classes and their centralizers in $\Sp(2n, \FF_q)$ due to Wall \cite{Wa62}.  An element in $\Sp(2n, \FF_q)$ which squares to the identity must have elementary divisors only of the form $x \pm 1$, and thus must be conjugate in $\GL(2n, \FF_q)$ to a diagonal element with only $1$'s and $-1$'s on the diagonal.  By \cite[p. 36, Case B]{Wa62}, the only such elements which are in $\Sp(2n, \FF_q)$ must have an even number $2k$ of eigenvalues equal to $1$, and so $2(n-k)$ eigenvalues equal to $-1$, and there is a unique such conjugacy class.  Also by \cite[p. 36, Case B]{Wa62}, an element in this conjugacy class has centralizer in $\Sp(2n, \FF_q)$ with size equal to $|\Sp(2k, \FF_q) \times \Sp(2(n-k), \FF_q)|$.  Using the fact that $|\Sp(2m, \FF_q)| = q^{m^2} \prod_{i=1}^m (q^{2i} - 1)$, we have that the size of the conjugacy class of $\Sp(2n, \FF_q)$ consisting of elements with $2k$ eigenvalues equal to $1$ and $2(n-k)$ eigenvalues equal to $-1$ and which square to the identity is
$$ \frac{|\Sp(2n, \FF_q)|}{|\Sp(2k, \FF_q) \times \Sp(2(n-k), \FF_q)|} = q^{2k(n-k)} \frac{\prod_{i=1}^n (q^{2i} - 1)}{\prod_{i=1}^k (q^{2i}-1) \prod_{i=1}^{n-k} (q^{2i}-1)} = q^{2k(n-k)} \binom{n}{k}_{q^2}.$$
Thus, the total number of elements in $\GSp(2n, \FF_q)$ such that $g^2 = I$ and $\mu(g) = 1$ is
\begin{equation} \label{FirstTerm}
\sum_{k=0}^n q^{2k(n-k)} \binom{n}{k}_{q^2}.
\end{equation}

Now suppose $g \in \GSp(2n, \FF_q)$ such that $g^2 = I$ and $\mu(g) = -1$.  In other words, $g^2 = -\mu(g) I$.  By \cite[Proposition 4]{Vi04}, there is a unique such conjugacy class in $\GSp(2n, \FF_q)$, and by \cite[Proposition 3.2]{Vi05}, the centralizer in $\GSp(2n, \FF_q)$ of any element in this conjugacy class has size $(q-1)|\GL(n, \FF_q)|$.  Therefore, the total number of elements in $\GSp(2n, \FF_q)$ such that $g^2 = I$ and $\mu(g)=-1$ is
\begin{equation} \label{SecondTerm}
\frac{|\GSp(2n, \FF_q)|}{(q-1)|\GL(n, \FF_q)|} = \frac{|\Sp(2n, \FF_q)|}{|\GL(n, \FF_q)|} = q^{n(n+1)/2} \prod_{i=1}^n (q^i + 1).
\end{equation}
From \cite[Theorem 3]{Go85} and \cite[Corollary 6.1]{Vi05}, the expression (\ref{SecondTerm}) is also the sum of the degrees of all irreducible characters of the symplectic group $\Sp(2n, \FF_q)$ when $q$ is odd.  Finally, the total number of elements in $\GSp(2n, \FF_q)$ which square to $I$ is the sum of (\ref{FirstTerm}) and (\ref{SecondTerm}), as desired.
\end{proof}

From Theorem \ref{GSpRealSum}, the sum of the degrees of the real-valued characters of $\GSp(2n, \FF_q)$ is a polynomial in $q$ of degree $n^2 + n$, while from \cite[Corollary 6.1]{Vi05}, the sum of the degrees of all of the irreducible characters of $\GSp(2n, \FF_q)$ is a polynomial in $q$ of degree $n^2 + n + 1$.  In this case, the ratio of the sum of all degrees to the sum of the degrees of real-valued characters is roughly $q$, which is quite different from the $\GL(n, \FF_q)$ case.

If we let $H = \Sp(2n, \FF_q)$, $q$ odd, then it follows from the Frobenius-Schur involution formula that the expression (\ref{FirstTerm}) is exactly
\begin{equation} \label{SpReal}
\sum_{k=0}^n q^{2k(n-k)} \binom{n}{k}_{q^2} = \sum_{\chi \in \Irr(H)} \ep(\chi) \chi(1) = \sum_{\chi \in \Irr(H) \atop{ \ep(\chi) = 1 }} \chi(1) - \sum_{\chi \in \Irr(H) \atop{ \ep(\chi) = -1}} \chi(1).
\end{equation}
In the case that $q \equiv 1($mod $4)$, it follows from a result of Wonenburger \cite[Theorem 2]{Wo66} and the fact that $-1$ is a square in $\FF_q$, that $\ep(\chi) = \pm 1$ for every irreducible character of $\Sp(2n, \FF_q)$.  When $q \equiv 3($mod $4)$, however, there are characters of $\Sp(2n, \FF_q)$ which are not real-valued, by \cite[Lemma 5.3]{FeZu82}.  Combining these facts with (\ref{SpReal}) and Theorem \ref{GSpRealSum}, we obtain the following. 

\begin{corollary} \label{SpSum}
Let $q$ be the power of an odd prime, let $H = \Sp(2n, \FF_q)$, and let $G = \GSp(2n, \FF_q)$.  Then,
$$ 2 \sum_{\chi \in \Irr(H) \atop{ \ep(\chi) = 1 }} \chi(1) + \sum_{\chi \in \Irr(H) \atop{ \ep(\chi) = 0 }} \chi(1) = \sum_{\chi \in {\rm Irr}(G) \atop{ \chi \; \RR \text{-valued}}} \chi(1).$$
If $q \equiv 1($mod $4)$, then for $\delta = 1$ or $-1$, we have
$$ \sum_{\chi \in \Irr(H) \atop{ \ep(\chi) = \delta }} \chi(1) = \frac{1}{2} \left(q^{n(n+1)/2} \prod_{i=1}^n (q^i + 1) + \delta \sum_{k=0}^n q^{2k(n-k)} \binom{n}{k}_{q^2} \right).$$
\end{corollary}

We now turn to the orthogonal similitude groups, $\GO^{\pm}(2n, \FF_q)$.  By Theorem \ref{GOReal}, every real-valued character of $\GO^{\pm}(2n, \FF_q)$ has Frobenius-Schur indicator $1$, and so again we find the sum of the degrees of the real-valued characters of this group by counting the number of elements which square to the identity.  Just as in the case of $\GSp(2n, \FF_q)$, such an element must either be orthogonal, or skew-orthogonal (that is, has similitude $-1$).  We prove that in the group $\GO^+(2n, \FF_q)$, the latter elements form a single conjugacy class, whereas there are no such elements in $\GO^-(2n, \FF_q)$.  These facts follow from the classification of conjugacy classes in $\GO^{\pm}(2n, \FF_q)$ due to Shinoda \cite{Sh80}, although we are also able to give an elementary proof, which is essentially the same as the proof of a result of Gow \cite[Lemma 1]{Go88} for symplectic groups.  We also find the centralizer of a skew-orthogonal order $2$ element in the group $\GO^+(2n, \FF_q)$.  The following result could certainly be obtained for a more general type of conjugacy class in the group of orthogonal similitudes over an arbitrary field (with characteristic not $2$), similar to the results for groups of symplectic similitudes in \cite[Proposition 4]{Vi04} and \cite[Proposition 3.2]{Vi05}, but we do not need such a result here.

\begin{lemma} \label{GOConjClass}  Let $G = \GO^{\pm}(2n, \FF_q)$, where $q$ is the power of an odd prime.  Consider the set $K = \{ g \in G \, \mid \, g^2 = I, \mu(g) = -1 \}$.  Then we have the following:
\begin{enumerate}
\item  If $G = \GO^+(2n, \FF_q)$, then the set $K$ forms a single conjugacy class in $G$.  The centralizer in $G$ of any element of $K$ is isomorphic to $\FF_q^{\times} \times \GL(n, \FF_q)$.
\item  If $G = \GO^-(2n, \FF_q)$, then the set $K$ is empty.
\end{enumerate}
\end{lemma}
\begin{proof}  Let $V = \FF_q^{2n}$, let $\langle \cdot, \cdot \rangle$ denote a nondegenerate symmetric form on $V$, and let $G$ be the orthogonal group of similitudes on $\langle \cdot, \cdot \rangle$.  If $g \in G$ such that $g^2 = I$ and $\mu(g) = -1$, then the only eigenvalues of $g$ are $1$ and $-1$.  Let $V_1$ and $V_{-1}$ be the eigenspaces of $g$ for $1$ and $-1$, respectively, and note that the assumption $\mu(g) = -1$ forces $V_1$ and $V_{-1}$ to be totally isotropic spaces with respect to $\langle \cdot, \cdot \rangle$.  This implies both $V_1$ and $V_{-1}$ have dimension $n$, and through the inner product $\langle \cdot, \cdot \rangle$, $V_{-1}$ is isomorphic to the dual space of $V_1$.  If $v_1, \ldots v_n$ is any basis of $V_1$, choose $w_1, \ldots, w_n$ to be a dual basis of $V_{-1}$, so that $\langle v_i, w_j \rangle = \delta_{ij}$.  With respect to this basis, the form $\langle \cdot, \cdot \rangle$ can be represented by the matrix $\left( \begin{array} {cc} 0 & I \\ I & 0 \end{array} \right)$, which means that it must be a split form on $V$.  That is, assuming that an element $g \in G$ with the above properties exists, we cannot have that $G$ is a group corresponding to a non-split form, which proves statement $2$.  We may now assume $G = \GO^+(2n, \FF_q)$.  With respect to the basis we have chosen, $g$ must be the element $\left( \begin{array} {cc} I & 0 \\ 0 & -I \end{array} \right)$, thus determining the conjugacy class of $g$ uniquely.

Now, without loss of generality, we may assume the form $\langle \cdot, \cdot \rangle$ is given by the matrix $\left( \begin{array} {cc} 0 & I \\ I & 0 \end{array} \right)$, and the element $g = \left( \begin{array} {cc} I & 0 \\ 0 & -I \end{array} \right)$.  It is a direct computation that the centralizer of $g$ in $G$ consists of the following set of elements:
$$ \left\{ \left( \begin{array} {cc} A & 0 \\ 0 & \lambda ({^T A}^{-1}) \end{array} \right) \, \mid \, A \in \GL(n, \FF_q), \lambda \in \FF_q^{\times} \right\}.$$
By mapping the element $\left( \begin{array} {cc} A & 0 \\ 0 & \lambda ({^T A}^{-1}) \end{array} \right)$ to $(\lambda, A)$, we see that this group is isomorphic to $\FF_q^{\times} \times \GL(n, \FF_q)$, completing the proof of statement 1.
\end{proof}

We now use the above result to give expressions for the sums of the degrees of real-valued characters of the groups $\GO^{\pm}(2n, \FF_q)$.  The key result due to Gow \cite[Theorem 1]{Go85} we use here, which we also used in the previous section, is that any irreducible character $\chi$ of any finite orthogonal group $\gO(n, \FF_q)$ (where $q$ is odd) satisfies $\ep(\chi) = 1$.

\begin{theorem} \label{GORealSum}  The sum of the degrees of the irreducible real-valued characters of the groups $\GO^{\pm}(2n, \FF_q)$, where $q$ is the power of an odd prime, are given as follows:
\begin{enumerate}
\item  If $G = \GO^+(2n, \FF_q)$ and $H = \gO^+(2n, \FF_q)$, then
$$\sum_{\chi \in {\rm Irr}(G) \atop{ \chi \; \RR \text{-valued}}} \chi(1)  = \sum_{ \chi \in {\rm Irr}(H)} \chi(1) + \frac{|\gO^+(2n, \FF_q)|}{|\GL(n, \FF_q)|}  = \sum_{ \chi \in {\rm Irr}(H)} \chi(1) + 2q^{n(n-1)/2} \prod_{i=1}^{n-1} (q^i + 1).$$
\item  If $G = \GO^-(2n, \FF_q)$, and $H = \gO^-(2n, \FF_q)$, then
$$ \sum_{\chi \in {\rm Irr}(G) \atop{ \chi \; \RR \text{-valued}}} \chi(1)  = \sum_{ \chi \in {\rm Irr}(H)} \chi(1).$$
\end{enumerate}
\end{theorem}
\begin{proof}  From Theorem \ref{GOReal}, the sum of the degrees of the real-valued characters of $\GO^{\pm}(2n, \FF_q)$ is equal to the number of elements $g$ in the group such that $g^2 = I$.  In the case $G = \GO^+(2n, \FF_q)$, we add the number of such elements in the group $\gO^+(2n, \FF_q)$ to the number of these elements such that $\mu(g) = -1$.  The number of elements in $\gO^+(2n, \FF_q)$ which square to the identity is equal to the sum of the degrees of all of the irreducible characters of $\gO^{+}(2n, \FF_q)$, by the Frobenius-Schur formula and \cite[Theorem 1]{Go85}.  By Lemma \ref{GOConjClass}(1), the number of order $2$ elements in $\GO^+(2n, \FF_q)$ such that $\mu(g) = -1$ is equal to the index of the centralizer of the unique conjugacy class of such elements in the group, which is given by
$$ \frac{|\GO^+(2n, \FF_q)|}{|\FF_q^{\times} \times \GL(n, \FF_q)|} = \frac{|\gO^+(2n, \FF_q)|}{|\GL(n, \FF_q)|} = \frac{ 2q^{n(n-1)} (q^n - 1) \prod_{i=1}^{n-1} (q^{2i} - 1)}{ q^{n(n-1)/2} \prod_{i=1}^n(q^i - 1)} = 2q^{n(n-1)/2} \prod_{i=1}^{n-1} (q^i + 1).$$
Adding this quantity to the sum of the degrees of the irreducible characters of $\gO^+(2n, \FF_q)$ gives the result for this case.

In the case $G = \GO^-(2n, \FF_q)$, by Lemma \ref{GOConjClass}(2), the only elements in $G$ which square to the identity are in the group $H = \gO^-(2n, \FF_q)$.  So, the sum of the degrees of the real-valued irreducible characters of $G$ is equal to the number of elements in $H$ which square to the identity, which, as in the previous case, is equal to the sum of the degrees of the irreducible characters of $H$.
\end{proof}

In Theorem \ref{OrthSum} below, we find expressions in $q$ for the sums of the degrees of the orthogonal groups $\gO(n, \FF_q)$, which, combined with Theorem \ref{GORealSum}, give us expressions for the sums of the degrees of the real-valued irreducible characters of $\GO^{\pm}(2n, \FF_q)$ as polynomials in $q$.

There are several similarities between Theorems \ref{GSpRealSum} and \ref{GORealSum}(1).  Both sums have a part corresponding to the sum of the degrees of the irreducible characters of the symplectic or orthogonal subgroups, respectively.  Also, both sums have a part which is the index of a general linear group as a subgroup of the symplectic or orthogonal subgroups, which corresponds to the size of a unique conjugacy class.  This part of the sum is conveniently factorizable as a polynomial in $q$.  The main difference is that this part of the sum corresponds to the sum of the degrees of the irreducible characters of the symplectic group in Theorem \ref{GSpRealSum}, but it does not correspond to the sum of the degrees for the orthogonal group in Theorem \ref{GORealSum}(1).  This transposed difference of the role of this part of the sum is essentially due to the difference between the defining bilinear forms, one being skew-symmetric, and the other symmetric, which switches the roles of the two different types of order $2$ elements which must be counted in each case.  As a result, we will see in the next section that the sum of the degrees of the irreducible characters of a finite orthogonal group is not, in general, a conveniently factorizable polynomial in $q$, but rather another sum involving Gaussian binomial coefficients. 

\section{Character degree sums for orthogonal groups} \label{OrthogSums}

In this section we compute an expression for the sum of the degrees of the irreducible characters of the orthogonal groups over finite fields of odd characteristic.  By the main theorem of \cite{Go85}, every irreducible character of such a group is the character of a real representation, and so by the classical Frobenius-Schur involution formula, the sum of the degrees of the irreducible characters is equal to the number of elements which square to the identity in the group.

So, in order to find the character degree sum for these finite orthogonal groups, we must count the number of involutions, which we do using the results on conjugacy classes in these groups due to Wall \cite{Wa62}.

If $g \in {\rm O}(n, \FF_q)$ and $g^2 = I$, then any elementary divisor of $g$ must be either $x+1$ or $x-1$, so over $\GL(n,\FF_q)$ is conjugate to a diagonal matrix with only $1$'s and $-1$'s on the diagonal.  By the results of Wall \cite[p. 36, Case B]{Wa62} (see also \cite[Section 2.2]{Fu00}), if there are $j$ eigenvalues of $g$ equal to $1$, and $n-j$ eigenvalues equal to $-1$, where $0 < j < n$, then there are exactly $2$ conjugacy classes in ${\rm O}(n, \FF_q)$ of such elements which square to the identity.  Let us label these conjugacy classes by the notation $(j, n-j)^{\pm}$, where $j$ is the number of eigenvalues equal to $1$, $n-j$ the number of eigenvalues equal to $-1$, and the sign $\pm$ distinguishes the two corresponding conjugacy classes.  The following summarizes the results we use from \cite{Wa62} for these conjugacy classes and the order of the centralizers of elements in these conjugacy classes.

\begin{proposition} [Wall] \label{OrthCent}
Let $G$ be an orthogonal group over $\FF_q$ with $q$ the power of an odd prime.
\begin{enumerate}
\item If $G = {\rm O}(n, \FF_q)$ with $n$ odd, and $j$ is even, then an element in the conjugacy class $(j, n-j)^{\pm}$ has centralizer size $|\gO^{\pm}(k, \FF_q) \times \gO(n-j, \FF_q)|$.  If $j$ is odd, an element in the conjugacy class $(j, n-j)^{\pm}$ has centralizer size $|\gO(j, \FF_q) \times \gO^{\pm}(n-j, \FF_q)|$.

\item If $G = \gO^+(n, \FF_q)$ with $n$ even, and $j$ is even, then an element in the conjugacy class $(j, n-j)^{\pm}$ has centralizer size $|\gO^{\pm}(j, \FF_q) \times \gO^{\pm}(n-j, \FF_q)|$.  If $j$ is odd, an element in either of the conjugacy classes $(j, n-j)^{\pm}$ has centralizer size $|\gO(j,\FF_q) \times \gO(n-j, \FF_q)|$.

\item If $G = \gO^-(n, \FF_q)$ with $n$ even, and $j$ is even, then an element in the conjugacy class $(j, n-j)^{\pm}$ has centralizer size $|\gO^{\pm}(j, \FF_q) \times \gO^{\mp}(n-j, \FF_q)|$.  If $j$ is odd, an element in either of the conjugacy classes $(j, n-j)^{\pm}$ has centralizer size $|\gO(j,\FF_q) \times \gO(n-j, \FF_q)|$.
\end{enumerate}
\end{proposition}

Using Proposition \ref{OrthCent}, and that the orders of the orthogonal groups are given by
$$ |\gO(2m+1, \FF_q)| = 2q^{m^2} \prod_{i=1}^m (q^{2i}-1), \quad \text{ and } \quad |\gO^{\pm}(2m, \FF_q)| = 2q^{m(m-1)} (q^m \mp 1) \prod_{i=1}^{m-1} (q^{2i} -1),$$
when $q$ is the power of an odd prime, we may obtain the following.

\begin{theorem} \label{OrthSum}
The sums of the character degrees of the orthogonal groups over $\FF_q$, where $q$ is the power of an odd prime, are given as follows:
\begin{enumerate}
\item If $G = {\rm O}(n, \FF_q)$, where $n = 2m+1$ is odd, then
$$ \sum_{\chi \in \Irr(G)} \chi(1) = 2 \sum_{k=0}^m q^{2k(m-k+1)} \binom{m}{k}_{q^2}.$$

\item If $G = {\rm O}^{\pm}(n, \FF_q)$, where $n=2m$ is even, then
$$ \sum_{\chi \in \Irr(G)} \chi(1) = \sum_{k=0}^m q^{2k(m-k)} \binom{m}{k}_{q^2} + q^{m-1} (q^m \mp 1) \sum_{k=0}^{m-1} q^{2k(m-k-1)} \binom{m-1}{k}_{q^2}.$$
\end{enumerate}
\end{theorem}
\begin{proof} In each case, we must sum the indices of the centralizers of the conjugacy classes of elements which square to the identity, using Proposition \ref{OrthCent}.  In case 1, when $G = \gO(2m+1, \FF_q)$, first consider the two conjugacy classes of type $(j, 2m+1-j)^{\pm}$, when $j=2k$ is even, $0 < k \leq m$.  Then, the sizes of the two centralizers $C^{\pm}$ of these classes are
$$|C^{\pm}| = |\gO^{\pm}(2k, \FF_q) \times \gO(2(m-k)+1, \FF_q)| = 4q^{k(k-1) + (m-k)^2} (q^k \mp 1) \prod_{i=1}^{k-1} (q^{2i}-1) \prod_{i=1}^{m-k} (q^{2i} - 1).$$
Computing the sum of the indices of these centralizers, we obtain
$$ \frac{|G|}{|C^{+}|} + \frac{|G|}{|C^{-}|} = q^{2k(m-k+1)} \frac{\prod_{i=1}^m (q^{2i} - 1)}{\prod_{i=1}^k (q^{2i} - 1) \prod_{i=1}^{m-k} (q^{2i}-1)} = q^{2k(m-k+1)} \binom{m}{k}_{q^2}.$$
Consider now the two conjugacy classes of type $(j, 2m+1-j)^{\pm}$, when $j = 2k+1$ is odd, so $2m+1 - j = 2(m-k)$, and $0 \leq k < m$.  The orders of the two centralizers $C^{\pm}$ in this case are
$$ |C^{\pm}| = |\gO(2k+1, \FF_q) \times \gO^{\pm}(2(m-k), \FF_q)|,$$
which are the same as in the previous case, except that $k$ is replaced by $m-k$.  Thus, the sum of the indices of these centralizers is
$$ \frac{|G|}{|C^+|} + \frac{|G|}{|C^-|} = q^{2(m-k)(k+1)} \binom{m}{m-k}_{q^2}.$$
Taking the sum of these terms, and replacing $k$ by $m-k$ in the sum, we obtain
$$ \sum_{k=0}^{m-1} q^{2(m-k)(k+1)} \binom{m}{m-k}_{q^2} = \sum_{k=1}^m q^{2k(m-k+1)} \binom{m}{k}_{q^2}$$
elements from these conjugacy classes which square to the identity.  Adding in the two central elements $\pm I$, we obtain that the total number of elements in $\gO(2m+1, \FF_q)$ which square to the identity, and so the sum of the degrees of the irreducible characters is
$$ 2 \sum_{k=1}^m q^{2k(m-k+1)} \binom{m}{k}_{q^2} + 2 = 2\sum_{k=0}^m q^{2k(m-k+1)} \binom{m}{k}_{q^2}.$$

Now consider case 2, when $G = \gO^{\pm}(2m, \FF_q)$.  If $G = \gO^+(2m, \FF_q)$, then the two conjugacy classes of type $(j, 2m-j)^{\pm}$, where $j = 2k$ is even, $0 < k < m$, have centralizers $C^{\pm}$ with orders $|\gO^{\pm}(2k, \FF_q) \times \gO^{\pm}(2(m-k), \FF_q)|$, so
$$|C^{\pm}| =  4 q^{k(k-1) + (m-k)(m-k-1)} (q^k \mp 1)(q^{m-k} \mp 1) \prod_{i=1}^{k-1}(q^{2i}-1) \prod_{i=1}^{m-k-1}(q^{2i}-1).$$ 
If $G = \gO^-(2m, \FF_q)$, the two corresponding conjugacy classes have centralizers $C^{\pm}$ with orders $|\gO^{\pm}(2k, \FF_q) \times \gO^{\mp}(2(m-k), \FF_q)|$, which is the same as above, except that $q^{m-k} \mp 1$ is replaced by $q^{m-k} \pm 1$.  In both cases, the sum of the indices of these centralizers is computed to be
$$ \frac{|G|}{|C^+|} + \frac{|G|}{|C^-|} = q^{2k(m-k)} \binom{m}{k}_{q^2}.$$
Adding in the two central elements, these conjugacy classes contribute exactly
\begin{equation} \label{EvenCon}
\sum_{k=1}^{m-1} q^{2k(m-k)} \binom{m}{k}_{q^2} + 2 = \sum_{k=0}^m q^{2k(m-k)} \binom{m}{k}_{q^2}.
\end{equation}
When $G = \gO^{\pm}(2m, \FF_q)$, the two conjugacy classes of type $(2k+1, 2m-(2k+1))^{\pm}$, $0 \leq k \leq m-1$, have centralizers of the same order $|\gO(2k+1, \FF_q) \times \gO(2(m-k-1) + 1, \FF_q)|$.  The union of these two conjugacy classes thus has cardinality
$$ \frac{2|\gO^{\pm}(2m, \FF_q)|}{|\gO(2k+1, \FF_q) \times \gO(2(m-k-1) + 1, \FF_q)|} = q^{2k(m-k-1) + m -1} (q^m \mp 1) \binom{m-1}{k}_{q^2}.$$
Taking the sum of these contributions from $k=0$ to $k=m-1$, and adding to (\ref{EvenCon}), the result is obtained. 
\end{proof}  

In the case that $n=2m+1$ is odd, the orthogonal group $\gO(n, \FF_q)$ is just the direct product $\{\pm I\} \times \SO(2m+1, \FF_q)$ of the center with the special orthogonal group.  The sum of the degrees of the characters of $\SO(2m+1, \FF_q)$ is thus exactly half the sum for $\gO(2m+1, \FF_q)$ obtained in case 1 of Theorem \ref{OrthSum}.  That is, we have the following.

\begin{corollary} \label{SOcor}
Let $q$ be the power of an odd prime, and let $G = \SO(2m+1, \FF_q)$.  Then
$$ \sum_{\chi \in \Irr(G)} \chi(1) = \sum_{k = 0}^m q^{2k(m-k+1)} \binom{m}{k}_{q^2}.$$
\end{corollary}

When $n = 2m$ is even, Gow \cite[Theorem 2]{Go85} proved that every real-valued character of $\SO^{\pm}(2m, \FF_q)$ is the character of a real representation, and in the case $n = 4l$ is divisible by $4$, every character of $\SO^{\pm}(4l, \FF_q)$ is real-valued.  To find the sum of the degrees of the real-valued characters of $\SO^{\pm}(2m, \FF_q)$ in these cases, then, we have to count the number of elements of the group which square to the identity.  This is exactly the number of elements in $\gO^{\pm}(2m, \FF_q)$ which square to the identity and which have determinant $1$, and so have an even number of eigenvalues equal to $-1$ (and an even number equal to $1$).  This is just the first part of the sum obtained in case 2 of Theorem \ref{OrthSum}.  Curiously, this is exactly what is obtained as the sum of the degrees of the real-valued characters of $\GL(m, \FF_{q^2})$ in (\ref{GLReal}).  We summarize these observations below.

\begin{corollary} \label{SOcor2}
Let $q$ be the power of an odd prime, let $H = \SO^{\pm}(2m, \FF_q)$, and let $G = \GL(m, \FF_{q^2})$.  Then,
$$\sum_{\chi \in {\rm Irr}(H) \atop{ \chi \; \RR \text{-valued}}} \chi(1) = \sum_{k=0}^m q^{2k(m-k)} \binom{m}{k}_{q^2} = \sum_{\chi \in {\rm Irr}(G) \atop{ \chi \; \RR \text{-valued}}} \chi(1).$$
In the case that $m = 2l$ is also even, this is the sum of the degrees of all of the characters of the group $\SO^{\pm}(4l, \FF_q)$.
\end{corollary}

\section{Inequality Lemmas} \label{IneqLemmas}

In this section, we prove several inequalities in preparation for the results in Section \ref{Bound}.  We begin with an elementary bound for the Gaussian binomial coefficients.

\begin{lemma} \label{binomineq} For any integers $m \geq 1$, $1 \leq k \leq m$, we have
$$ \binom{m}{k}_q \leq q^{k(m-k) - m+1} (q+1)^{m-1}.$$
\end{lemma} 
\begin{proof}  The inequality reduces to $1 \leq 1$ when $m = 1$.  If $k = 1$, then for any $m \geq 1$, we have $\binom{m}{1}_q = q^{m-1} + \cdots + q + 1 \leq (q+1)^{m-1}$.  So, we may assume $k \geq 2$.  Assume the inequality holds for $m=n-1$, for any $k \geq 1$.  We use the identity $\binom{n}{k}_q = \binom{n-1}{k}_q + q^{n-k}\binom{n-1}{k-1}_q$, $k \geq 1$, which was mentioned in Section \ref{realsums}.  By this identity and the induction hypothesis, and with $k \geq 2$, we have
\begin{align*}
\binom{n}{k}_q & \leq (q^{k(n-1-k)-n+2} + q^{k(n-k) -n+2})(q+1)^{n-2} = q^{k(n-1-k) - n+2} (q^k + 1)(q+1)^{n-2} \\
               & \leq q^{k(n-1-k)-n+2}q^{k-1}(q+1)(q+1)^{n-2} = q^{k(n-k) - n+1} (q+1)^{n-1},
\end{align*}
as desired.
\end{proof}

The next two Lemmas will be used to bound the expressions obtained in Section \ref{OrthogSums}.

\begin{lemma} \label{EvenDimIneq}
For any integer $m \geq 1$, and any $q > 1$,
$$ \sum_{k=0}^m q^{2k(m-k)} \binom{m}{k}_{q^2} \leq \left\{ \begin{array}{ll} 2(q+1)^{m^2-1} & \text{if $m$ is odd} \\ (q+1)^{m^2} & \text{if $m$ is even.} \end{array} \right.$$
\end{lemma}
\begin{proof}
From the symmetry in $k$ and $m-k$ in the sum, we have
\begin{equation} \label{SymmSum}
\sum_{k=0}^m q^{2k(m-k)} \binom{m}{k}_{q^2} = \left\{ \begin{array}{ll} 2\sum_{k=0}^{(m-1)/2} q^{2k(m-k)} \binom{m}{k}_{q^2} & \text{if $m$ is odd} \\ q^{m^2/2} \binom{m}{m/2}_{q^2} + 2\sum_{k=0}^{(m/2) - 1} q^{2k(m-k)} \binom{m}{k}_{q^2} & \text{if $m$ is even.} \end{array} \right.
\end{equation}
From Lemma \ref{binomineq}, we have, for any $k \geq 1$, 
$$ \binom{m}{k}_{q^2} \leq q^{2k(m-k)-2m+2} (q^2 + 1)^{m-1} \leq q^{2k(m-k)-m+1} (q+1)^{m-1},$$
since $q^2 + 1 \leq q(q+1)$.  Now, for $s = (m-1)/2$ if $m$ is odd, or $s = (m/2) - 1$ if $m$ is even, we have
\begin{align} \label{Intermed}
\sum_{k=0}^s q^{2k(m-k)} \binom{m}{k}_{q^2} & \leq (q+1)^{m-1} \left( 1 + \sum_{k=1}^s q^{4k(m-k) - m+1} \right) \notag\\
& \leq (q+1)^{m-1} (q+1)^{4s(m-s) - m+1} = (q+1)^{4s(m-s)},
\end{align}
since the exponent of $q$ in the sum is maximum when $k=s$.  When $m$ is odd, then substituting $s=(m-1)/2$ and applying (\ref{SymmSum}) gives the desired result.  When $m$ is even, then applying (\ref{Intermed}) with $s = (m/2) - 1$, and Lemma \ref{binomineq}, we have
\begin{align*}
q^{m^2/4} \binom{m}{m/2}_{q^2} + 2\sum_{k=0}^{(m/2) - 1} q^{2k(m-k)} \binom{m}{k}_{q^2}  & \leq q^{m^2 - m + 1}(q+1)^{m-1} + 2(q+1)^{m^2 - 4} \\
    & \leq (q+1)^{m^2 - 4}(q^4 + 2) \leq (q+1)^{m^2}, 
\end{align*}
as claimed for $m$ even.
\end{proof}

\begin{lemma} \label{OddDimIneq}
For any integer $m \geq 0$, and any $q > 1$,
$$ \sum_{k=0}^m q^{2k(m-k+1)} \binom{m}{k}_{q^2} \leq (q+1)^{m^2 + m}.$$
\end{lemma}
\begin{proof} We may assume $m \geq 1$.  Note that by switching the roles of $k$ and $m-k$, we obtain
$$ \sum_{k=0}^m q^{2k(m-k+1)} \binom{m}{k}_{q^2} = \sum_{k=0}^m q^{2k} q^{2k(m-k)} \binom{m}{k}_{q^2} = \sum_{k=0}^m q^{2(m-k)} q^{2k(m-k)} \binom{m}{k}_{q^2}.$$
It follows that we have
$$2\sum_{k=0}^m q^{2k(m-k+1)} \binom{m}{k}_{q^2} = \sum_{k=0}^m (q^{2k}+q^{2(m-k)})q^{2k(m-k)} \binom{m}{k}_{q^2}.$$
From the symmetry in $k$ and $m-k$ in the right-hand side of the above equation, we have
\begin{equation*}
\sum_{k=0}^m q^{2k(m-k+1)} \binom{m}{k}_{q^2} = 
\end{equation*}
\begin{equation} \label{SymmSum2}
=\left\{ \begin{array}{ll} \sum_{k=0}^{(m-1)/2}(q^{2k} + q^{2(m-k)})q^{2k(m-k)}\binom{m}{k}_{q^2} & \text{if $m$ is odd} \\ q^{m^2/2 + m}\binom{m}{m/2}_{q^2} + \sum_{k=0}^{(m/2)-1} (q^{2k} + q^{2(m-k)}) q^{2k(m-k)} \binom{m}{k}_{q^2} & \text{if $m$ is even}. \end{array} \right.
\end{equation}
Let $s = (m-1)/2$ if $m$ is odd, and $s = (m/2)-1$ if $m$ is even.  In the sums (\ref{SymmSum2}), the term $q^{2k} + q^{2(m-k)}$ takes its maximum value when $k=s$.  Applying this, and the inequality (\ref{Intermed}), we have
\begin{align} \label{Intermed2}
\sum_{k=0}^s (q^{2k} + q^{2(m-k)}) q^{2k(m-k)} \binom{m}{k}_{q^2} & \leq (q^{2s} + q^{2(m-s)}) \sum_{k=0}^s q^{2k(m-k)} \binom{m}{k}_{q^2} \notag \\
& \leq (q^{2s} + q^{2(m-s)})(q+1)^{4s(m-s)}.
\end{align}
When $m$ is odd and $s = (m-1)/2$, then by (\ref{Intermed2}) and (\ref{SymmSum2}), 
$$\sum_{k=0}^m q^{2k(m-k+1)} \binom{m}{k}_{q^2} \leq (q^{m-1} + q^{m+1})(q+1)^{(m^2-1)/2} = q^{m-1}(q+1)^{m^2 + 1} \leq (q+1)^{m^2 + m},$$
as claimed.  When $m$ is even, with $s = (m/2) -1$, apply (\ref{Intermed2}), together with (\ref{SymmSum2}) and Lemma \ref{binomineq} to obtain
\begin{align*}
\sum_{k=0}^m q^{2k(m-k+1)} \binom{m}{k}_{q^2} & \leq q^{m^2 + 1} (q+1)^{m-1} + (q^{m-2} + q^{m+2})(q+1)^{m^2 - 4} \\
& \leq q^{m^2 + 1} (q+1)^{m-1} + q^{m-2} (q^4 + 1)(q + 1)^{m^2 - 4} \\
& \leq q^2 (q+1)^{m^2 + m - 2} + (q+1)^{m^2 +m -2} \leq (q+1)^{m^2+m},
\end{align*}
as desired.
\end{proof}

\section{An upper bound for character degree sums} \label{Bound}

Let $\bG$ be a connected algebraic group with connected center over $\bar{\FF}_q$, defined over $\FF_q$ by some Frobenius map.  By the {\em dimension} of $\bG$, we mean the dimension of $\bG$ as an algebraic variety over $\bar{\FF}_q$.  The {\em rank} of $\bG$ is the dimension of a maximal torus of $\bG$.  That is, if the rank of $\bG$ is $r$, then a maximal torus of $\bG$ is isomorphic to $(\bar{\FF}_q^{\times})^r$.  

In this section, we will consider the case when $\bG$ is a connected {\em classical} group with connected center, when $q$ is the power of an odd prime.  For us, these are the groups $\GL(n, \bar{\FF}_q)$, $\SO(2n+1, \bar{\FF}_q)$, $\GSp(2n, \bar{\FF}_q)$, and a certain index $2$ subgroup of $\GO(2n, \bar{\FF}_q)$, denoted $\GO^{\circ}(2n, \bar{\FF}_q)$, which we describe now.

Let $V$ be a $2n$-dimensional vector space over $\bar{\FF}_q$, with $q$ odd, and consider a nondegenerate symmetric form on $V$.  Since $\bar{\FF}_q$ is algebraically closed, there is only one equivalence class of forms (by \cite[Theorem 4.4]{Gr02}, for example), which we may assume corresponds to the matrix $J = \left( \begin{array} {cc} 0 & I \\ I & 0 \end{array} \right)$.  Consider the orthogonal group of similitudes $\GO(2n, \bar{\FF}_q)$ of this form.  For any $g \in \GO(2n, \bar{\FF}_q)$, then, we have ${^T g} = J \mu(g) g^{-1} J$.  Since $g$ is conjugate to its transpose in $\GL(V)$, it follows that $g$ is conjugate to $\mu(g) g^{-1}$ in $\GL(V)$, and so they have the same determinant.  In particular, we have ${\rm det}(g)^2 = \mu(g)^{2n}$.  The subgroup consisting of elements with the property that ${\rm det}(g) = \mu(g)^n$ is a connected algebraic group (by an argument similar to that given in \cite[Section 15.2]{DiMi91}).  So, we define
\begin{equation} \label{GOconnDefn}
\GO^{\circ} (2n, \bar{\FF}_q) = \{ g \in \GO(2n, \bar{\FF}_q) \, \mid \, {\rm det}(g) = \mu(g)^n \},
\end{equation}
the connected component of the identity in the orthogonal group of similitudes.  Note that $\SO(2n, \bar{\FF}_q)$ is contained in $\GO^{\circ}(2n, \bar{\FF}_q)$, and its center consists of all scalar matrices.

If we define the $\FF_q$-structure of $\GO^{\circ}(2n, \bar{\FF}_q)$ by the standard Frobenius map $F$, which raises entries to the power of $q$, we get that the group of $\FF_q$-points is the index $2$ subgroup of the split orthogonal group of similitudes over $\FF_q$ satisfying the condition in (\ref{GOconnDefn}), which we denote $\GO^{+, \circ}(2n, \FF_q)$.  If we compose the standard Frobenius map with conjugation by an orthogonal reflection defined over $\FF_q$ (see \cite[Section 15.3]{DiMi91}), we get that the group of $\FF_q$-points of $\GO^{\circ}(2n, \bar{\FF}_q)$ is the index $2$ subgroup of the non-split orthogonal group of similitudes over $\FF_q$ which satisfies (\ref{GOconnDefn}), which we denote $\GO^{-, \circ}(2n, \FF_q)$.  We may relate the character degree sums for the groups $\GO^{\pm, \circ}(2n, \FF_q)$ to those for the orthogonal groups $\gO^{\pm}(2n, \FF_q)$ in the following way.

\begin{lemma} \label{GOSumLemma} Let $q$ be the power of an odd prime, let $G = \GO^{\pm, \circ}(2m, \FF_q)$, and let $H = \gO^{\pm}(2m, \FF_q)$.  Then
$$ \sum_{ \chi \in {\rm Irr}(G)} \chi(1) \leq (q-1) \sum_{\chi \in {\rm Irr}(H)} \chi(1).$$
\end{lemma}
\begin{proof} Let $S = \SO^{\pm}(2m, \FF_q)$ be the special orthogonal group.  Then $S$ is a normal subgroup of $G$ of index $q-1$, with $G/S \cong \FF_q^{\times}$ cyclic.  From Clifford theory (see \cite[Chapter 6 and Theorem 11.7]{Is76}), the restriction of any irreducible character of $G$ to $S$ has a multiplicity-free decomposition, and each irreducible character of $S$ appears in the restriction to $S$ of at most $q-1$ characters of $G$.  It follows that we have
\begin{equation} \label{IneqOne}
\sum_{\chi \in {\rm Irr}(G)} \chi(1) \leq (q-1) \sum_{\chi \in {\rm Irr}(S)} \chi(1).
\end{equation}
Since $S$ is an index $2$ subgroup of $H$, every irreducible character of $H$ is extended or induced from an irreducible character of $S$.  A character of $H$ which is extended from a character of $S$ has the same degree as that character of $H$, while that character of $S$ may be extended to give a second distinct irreducible of $H$.  A character of $H$ which is induced from a character of $S$ has degree twice that of the character of $S$, but can also be induced by exactly one other distinct character of $S$.  It follows that we have
\begin{equation} \label{IneqTwo}
\sum_{\chi \in {\rm Irr}(S)} \chi(1) \leq \sum_{\chi \in {\rm Irr}(H)} \chi(1).
\end{equation}
The result follows from (\ref{IneqOne}) and (\ref{IneqTwo}).
\end{proof}

We may now prove the main result of this section, in which we improve Kowalski's Theorem \ref{KoThm} in the case of any connected classical group with connected center, as follows.

\begin{theorem} \label{MainBound} Let $q$ be the power of an odd prime, and let $\bG$ be a connected classical group with connected center defined over $\FF_q$, where the rank of $\bG$ is $r$, and the dimension of $\bG$ is $d$.  Then the sum of the degrees of the irreducible characters of the finite group $\bG(\FF_q)$ may be bounded as follows:
$$ \sum_{\chi \in {\rm Irr}(\bG(\FF_q))} \chi(1) \leq (q+1)^{(d+r)/2}.$$
\end{theorem}
\begin{proof} The dimensions and ranks of finite classical groups may be computed directly from their definitions (see \cite[Chapter 15] {DiMi91}, for example).  If $\bG = \GL(n, \bar{\FF}_q)$, then $d = n^2$ and $r = n$.  In this case, $\bG$ may have $\FF_q$-structure given by the standard Frobenius map, in which case $\bG(\FF_q) = \GL(n, \FF_q)$, or the standard Frobenius map composed with the transpose-inverse automorphism, in which case $\bG(\FF_q) = {\rm U}(n, \FF_{q^2})$, the finite unitary group defined over $\FF_q$.  The case $\bG(\FF_q) = \GL(n, \FF_q)$ was considered by Kowalski \cite[p. 80]{Ko08}, and he showed that the sum of the degrees of the characters is indeed bounded above by $(q+1)^{(d+r)/2}$.  When $\bG(\FF_q) = {\rm U}(n, \FF_{q^2})$, the sum of the degrees of the irreducible characters of $\bG(\FF_q)$ was computed by Thiem and the author \cite[Theorem 5.2]{ThVi07} to be
\begin{align*}
\sum_{\chi \in {\rm Irr}(\bG(\FF_q))} \chi(1)  & = (q+1)q^2(q^3+1)q^4 \cdots (q^n + (1-(-1)^n)/2)\\
& \leq \prod_{i=1}^n (q+1)^i = (q+1)^{(n^2 + n)/2} = (q+1)^{(d+r)/2}.
\end{align*}

If $\bG = \GSp(2n, \bar{\FF}_q)$, then $d = 2n^2 + n +1$ and $r = n+1$, and we may assume $\bG$ has $\FF_q$-structure given by the standard Frobenius map, so that $\bG(\FF_q) = \GSp(2n, \FF_q)$.  This case was also considered by Kowalski, and he noticed that by the formula given in \cite[Corollary 6.1]{Vi05}, the sum of the degrees of $\GSp(2n, \FF_q)$ is bounded above by $(q+1)^{(d+r)/2}$.  We note that the results in \cite{Vi05} are proven only when $q$ is odd, and so we can only obtain this bound in the case $q$ is odd.

If $\bG = \SO(2n+1, \bar{\FF}_q)$, then $d = 2n^2 + n$ and $r = n$, and again we may assume $\bG$ has $\FF_q$-structure given by the standard Frobenius map, and so $\bG(\FF_q) = \SO(2n+1, \FF_q)$.  By Corollary \ref{SOcor} and Lemma \ref{OddDimIneq}, we have
$$ \sum_{\chi \in {\rm Irr}(\bG(\FF_q))} \chi(1) = \sum_{k=0}^n q^{2k(n-k+1)} \binom{n}{k}_{q^2} \leq (q+1)^{n^2 + n} = (q+1)^{(d+r)/2}. $$

If $\bG = \GO^{\circ}(2n, \bar{\FF}_q)$, then $d = 2n^2 - n + 1$ and $r = n+1$.  In this case, as explained above, $\bG$ can have $\FF_q$-structure given by the standard Frobenius map, so that $\bG(\FF_q) = \GO^{+, \circ}(2n, \FF_q)$, or by the standard Frobenius map composed with conjugation by an orthogonal reflection defined over $\FF_q$, in which case $\bG(\FF_q) = \GO^{-, \circ}(2n, \FF_q)$.  We consider both cases at once, so let $\bG(\FF_q) = \GO^{\pm, \circ}(2n, \FF_q)$.  By Theorem \ref{OrthSum}(2) and Lemma \ref{GOSumLemma}, we have
\begin{equation} \label{GOIntermed}
\sum_{\chi \in {\rm Irr}(\bG(\FF_q))} \chi(1) \leq (q-1) \left(\sum_{k=0}^n q^{2k(n-k)} \binom{n}{k}_{q^2} + q^{n-1} (q^n + 1) \sum_{k=0}^{n-1} q^{2k(n-k-1)} \binom{n-1}{k}_{q^2} \right).
\end{equation}
If $n$ is even, then by (\ref{GOIntermed}) and Lemma \ref{EvenDimIneq}, we have
\begin{align*}
\sum_{\chi \in {\rm Irr}(\bG(\FF_q))} \chi(1) & \leq (q-1)\left( (q+1)^{n^2} + q^{n-1}(q^n + 1) 2(q+1)^{(n-1)^2-1} \right) \\
 & \leq (q-1) \left( (q+1)^{n^2} + 2(q+1)^{n^2 - 1} \right) = (q+1)^{n^2 - 1}(q-1)(q+3) \\
 & = (q+1)^{n^2 -1} (q^2 + 2q -3) \leq (q+1)^{n^2 + 1} = (q+1)^{(d+r)/2}, 
\end{align*}
as required.  Similarly, when $n$ is odd, then by (\ref{GOIntermed}) and Lemma \ref{EvenDimIneq}, we have
\begin{align*}
\sum_{\chi \in {\rm Irr}(\bG(\FF_q))} \chi(1) & \leq (q-1)\left( 2(q+1)^{n^2-1} + q^{n-1}(q^n + 1)(q+1)^{(n-1)^2} \right) \\
 & \leq (q-1) \left( 2(q+1)^{n^2-1} + (q+1)^{n^2} \right) = (q+1)^{n^2 - 1}(q-1)(q+3) \\
 & = (q+1)^{n^2 -1} (q^2 + 2q -3) \leq (q+1)^{n^2 + 1} = (q+1)^{(d+r)/2}, 
\end{align*}
which was the last case to consider.
\end{proof}

\section{A lower bound and a conjecture}  \label{Lower}

We now consider the case that $\bG$ is any connected reductive group over $\bar{\FF}_q$ with connected center, defined over $\FF_q$, where $q$ is the power of some prime $p$ (we allow $p=2$ here).  If $G = \bG(\FF_q)$, and $N$ is a maximal unipotent subgroup of $G$, then $N$ is a Sylow $p$-subgroup of $G$ \cite[Proposition 3.19(i)]{DiMi91}.  The {\em Gelfand-Graev character} of $G$ is constructed by inducing a {\em nondegenerate} linear character from $N$ to $G$.  For a detailed discussion on Gelfand-Graev characters and nondegenerate characters, see \cite[Section 8.1]{Ca85} or \cite[Chapter 14]{DiMi91}.  The main result on the Gelfand-Graev character which we need is that its decomposition into a sum of irreducible characters of $G$ is multiplicity-free, a result which is due to Steinberg \cite{St67} in the general case.  By applying this fact, we obtain the following lower bound for the sum of the degrees of the irreducible characters of $G = \bG(\FF_q)$.

\begin{proposition} \label{LowerProp} Let $\bG$ be a connected reductive group over $\bar{\FF}_q$ with connected center, defined over $\FF_q$, of dimension $d$ and rank $r$.  Then the sum of the dimensions of the irreducible representations of the group $G = \bG(\FF_q)$ is bounded below as follows:
$$ \sum_{\chi \in {\rm Irr}(G)} \chi(1) \geq  q^{(d-r)/2} (q-1)^r.$$
\end{proposition}
\begin{proof} Since the Gelfand-Graev character of $G = \bG(\FF_q)$ is multiplicity-free, then the degree of the Gelfand-Graev character gives a lower bound for the sum of the degrees of all of the irreducible characters of $G$.  We give a lower bound for the degree of the Gelfand-Graev character to prove our claim.

The {\em semisimple rank} of $\bG$, which we denote by $l$, is defined as the rank of $\bG/R(\bG)$, where $R(\bG)$ is the radical of $\bG$, which is the maximal closed, connected, normal, solvable subgroup of $\bG$ (see \cite[6.4.14]{Sp98}).  The degree of the Gelfand-Graev character is the $p'$-part of the order of $G$, since it is induced from a linear character of a Sylow $p$-subgroup of $G$.  By \cite[Section 2.9]{Ca85}, the $p'$-part of the order of $G$ is given by
$$ |\bZ(\FF_q)| \prod_{i=1}^l (q^{d_i} - \omega_i),$$
where $\bZ$ is the center of $\bG$, the $d_i$ are the degrees of the generators of the ring of polynomial invariants of the Weyl group $W$ of $\bG$, and the $\omega_i$ are roots of unity which are eigenvalues of a linear map on the algebra of polynomial invariants of $W$.  Now, we have $|q^{d_i} - \omega_i| \geq q^{d_i} - 1 \geq q^{d_i - 1}(q-1)$.  From \cite[Proposition 2.4.1(iv)]{Ca85}, we have $\sum_{i=1}^l (d_i - 1)$ is equal to the number of positive roots of the root system corresponding to the Weyl group $W$, and by \cite[Corollary 8.1.3(ii)]{Sp98}, this is equal to $(d-r)/2$.  These observations give 
\begin{equation} \label{firstpart}
\prod_{i=1}^l (q^{d_i} - \omega_i) \geq q^{\sum_{i=1}^l (d_i - 1)} (q-1)^l = q^{(d-r)/2} (q-1)^l.
\end{equation}

By \cite[Proposition 7.3.1(i)]{Sp98}, since $\bG$ is assumed to be a connected reductive group with connected center $\bZ$, we have $R(\bG) = \bZ$, where $\bZ$ is itself a torus.  In particular, we have $r = l + {\rm dim}(\bZ)$, by the definition of semisimple rank.  By \cite[Proposition 3.3.7]{Ca85} and the formula for the order of the $\FF_q$-points of a maximal torus \cite[Proposition 3.3.5]{Ca85}, like in \cite[p. 75]{Ko08}, this implies we have $|\bZ(\FF_q)| \geq (q-1)^{r - l}$.  Combining this with (\ref{firstpart}) gives us the desired bound.
\end{proof}

Let us return to the question of a general upper bound for the sum of the degrees of the irreducible characters of $\bG(\FF_q)$.  When $\bG(\FF_q) = \GL(n, \FF_q)$, then the sum of the character degrees of $\bG(\FF_q)$ is given by (see \cite{Go83, Kl83, Mac95})
$$(q-1)q^2(q^3-1) \cdots (q^n - (1 - (-1)^n)/2),$$
which is, of course, bounded above by $q^{(n^2+n)/2} < q^{(n^2 - n)/2}(q+1)^n$.  When $\bG(\FF_q) = {\rm U}(n, \FF_{q^2})$, then by the sum of the character degrees given in the proof of Theorem \ref{MainBound}, this can be seen to be bounded above by $q^{(n^2 - n)/2}(q+1)^n$ as well.  Note that in these two cases, we have $(n^2 - n)/2 = (d-r)/2$ and $n = r$.

When $\bG(\FF_q) = {\rm GSp}(2n, \FF_q)$, with $q$ odd, the sum of the degrees of the irreducible characters of $\bG(\FF_q)$, by \cite[Corollary 6.1]{Vi05}, is
$$ \frac{1}{2} q^{(n^2 + n)/2}(q-1) \left( \prod_{i=1}^n (q^i + 1) + \prod_{i=1}^n (q^i + (-1)^i) \right) \leq q^{n^2} (q+1)^{n+1},$$
where $n+1=r$ and $n^2 = (d-r)/2$.

If $\bG(\FF_q)$ is $\SO(2n+1, \FF_q)$ or $\GO^{\pm, \circ}(2n, \FF_q)$, with $q$ odd, then we may check directly for small values of $n$ that the sum of the character degrees of these groups are bounded by $q^{(d-r)/2}(q+1)^r$ as well.  We could perhaps prove this for all $n$ for these groups, although we would have to tighten the bounds found in Section \ref{IneqLemmas}.  It would be much more satisfying to have a general proof using Deligne-Lusztig theory, along the same lines as Kowalski's proof of Theorem \ref{KoThm} in \cite{Ko08}.

Based on these examples, we conclude by making the following conjecture for an upper bound for the sum of the degrees of the irreducible characters of a finite reductive group.  We include the lower bound obtained in Proposition \ref{LowerProp} in the statement to stress the symmetry in these bounds.

\begin{conjecture}  Let $\bG$ be a connected reductive group over $\bar{\FF}_q$ with connected center, defined over $\FF_q$, of dimension $d$ and rank $r$.  Then the sum of the degrees of the irreducible complex characters of $G = \bG(\FF_q)$ may be bounded as follows:
$$ q^{(d-r)/2}(q-1)^r \leq \sum_{\chi \in {\rm Irr}(G)} \chi(1) \leq q^{(d-r)/2}(q+1)^r.$$
\end{conjecture}

\bigskip

\noindent
\begin{tabular}{ll}
\textsc{Mathematics Department}\\ 
\textsc{College of William and Mary}\\
\textsc{P. O. Box 8795}\\
\textsc{Williamsburg, VA  23187}\\
{\em e-mail}:  {\tt vinroot@math.wm.edu}\\
\end{tabular}

\end{document}